\documentclass[10pt]{article}
\usepackage[utf8]{inputenc}
\usepackage[T1]{fontenc}
\usepackage[top=1in,bottom=1in,left=1in,right=1in]{geometry}
\usepackage{authblk}
\usepackage{amsmath,amsfonts,amssymb,amsthm}
\usepackage{hyperref}

\usepackage{tikz-cd}

\newcommand{\stab}{\mathfrak{stab}}
\newcommand{\struc}{\mathfrak{stab}^*}
\newcommand{\ustruc}{\mathfrak{stab}^*}
\newcommand{\skr}{\mathfrak{sK}_r}

\newtheorem{proposition}{Proposition}

\newtheorem{remark}{Remark}
\numberwithin{equation}{section}

\newcommand{\R}{\mathbb{R}}

\newcommand{\bin}[2]{\left(\begin{matrix} #1 \\ #2 \end{matrix}\right)}
\newcommand{\pfq}[3]{\left(\begin{matrix} #1 \\ #2 \end{matrix}\Bigg\rvert#3 \right)}

\title{The rational Sklyanin algebra and the Wilson and para-Racah polynomials}
 \author{Geoffroy Bergeron, Julien Gaboriaud, Luc Vinet and Alexei Zhedanov}
\date{\today}
\begin{document} 
\maketitle
\thispagestyle{empty}
\hrule
\begin{abstract}
    The relation between Wilson and para-Racah polynomials and representations of the degenerate rational Sklyanin algebra is established. Second order Heun operators on quadratic grids with no diagonal terms are determined. These special or S--Heun operators lead to the rational degeneration of the Sklyanin algebra; they also entail the contiguity and structure operators of the Wilson polynomials. The finite-dimensional restriction yields a representation that acts on the para-Racah polynomials.
\end{abstract}
\hrule

\section{Introduction}
This paper pursues the exploration of the links between Heun operators, Sklyanin algebras and orthogonal polynomials. Originally introduced in the context of quantum integrable systems \cite{sklyanin1984}, Sklyanin algebras are typically presented in terms of generators verifying homogeneous quadratic relations. These algebras have been the object of much attention from the perspective of algebraic geometry \cite{smith1994,walton2012,iyudu2017}. Classes of Heun operators can be defined \cite{grunbaum2018} from the property that they increase by no more than one the degree of polynomials defined on certain continuous or discrete domains; they have been the focus of a continued research effort \cite{takemura2018, baseilhac2019, vinet2019, crampe2019, tsujimoto2019, baseilhac2020, bergeron2020a} with many applications \cite{slepian1983,landau1985,bergeron2020,crampe2019a,crampe2019b,belliard2015,baseilhac2019a}. A key observation for our purposes is that a special category of these operators, referred to as S--Heun operators, offers a path towards the identification of interesting Sklyanin-like algebras through the relations they realize. This connects with orthogonal polynomials as these concrete S--Heun operators are recognized as ladder and structure operators for families of bispectral polynomials belonging to the Askey scheme. It is thus observed that these sets of orthogonal polynomials form representation bases for Sklyanin algebras. Furthermore, the finite-dimensional representations of these Sklyanin algebras are found to provide the algebraic setting that had so far been lacking for the orthogonal polynomials of the so-called ``para'' type. 

A first illustration of these connections was achieved in \cite{gaboriaud2020}. Building on results of Gorsky and Zabrodin \cite{gorsky1993} on the one hand and of Kalnins and Miller \cite{kalnins1989} on the other, this paper focused on S--Heun operators attached to the Askey--Wilson grid. The salient observations were: \textit{i.} that a subset of the S--Heun operators realize the trigonometric degeneration of the original elliptic Sklyanin algebra and \textit{ii.} that this Sklyanin algebra is a basic structure underneath the theory of Askey--Wilson polynomials.  Indeed, as was stressed, the Askey--Wilson operator admits a factorization in terms of the S--Heun operators realizing this degenerate Sklyanin algebra and as was also pointed out, the ladder and structure operators for the Askey--Wilson polynomials obtained by  Kalnins and Miller actually realize this degenerate algebra. In view of the fact that the Askey--Wilson algebra \cite{zhedanov1991} accounts for the bispectrality of the eponym polynomials, a parallel was thus drawn with the dynamical extension of symmetry algebras by the inclusion of ladder operators in the set of generators. Finally, the $q$-para Racah polynomials were seen to form a basis for the finite-dimensional representation of the degenerate Sklyanin algebra. This set the course for the systematic examination of the Sklyanin-like operators formed by S--Heun operators on lattices admitting orthogonal polynomials.

The study of S--Heun operators on linear and exponential grid and of the Sklyanin algebras they realize was carried out in \cite{bergeron2021}. It allowed to tie the representations of these algebras to the continuous Hahn and big $q$-Jacobi polynomials and in finite dimensions to the para-Krawtchouk and $q$-para Krawtchouk polynomials. This analysis confirmed the important role that Sklyanin algebras play in the interpretation of hypergeometric orthogonal polynomials. 

We here address the connection that the Wilson polynomials have with Sklyanin algebras. (We recall that these polynomials are at the top of the $q=1$ part of the Askey scheme.) This will call for the determination of the S--Heun operators on quadratic grids. The rational degeneration of the Sklyanin algebra first found by Smirnov \cite{smirnov2010} will be seen to emerge and to be realized by the structure and ladder operators \cite{miller1987} of the Wilson polynomials. This will hence attach these polynomials to representations of the rational Sklyanin algebra. In keeping with preceding observations, the finite-dimensional restrictions of these representations will be seen to offer an algebraic interpretation of the para-Racah polynomials \cite{lemay2016}.

\subsection{The Wilson polynomials and its truncations}
As the Wilson polynomials will prove central in deriving subsequent results, some of their known properties are summarized here. The four-parameter Wilson polynomials \cite{koekoek2010} of degree $n$, denoted $W_n(x^2\rvert a,b,c,d)$, are given by
\begin{align*}
 W_n(x^2\rvert a,b,c,d) = (a+b)_n(a+c)_n(a+d)_n\, {}_4F_3\pfq{-n,n+a+b+c+d-1,a+ix,a-ix}{a+b,a+c,a+d}{1},
\end{align*}
where $(a)_n=a(a+1)...(a+n-1)$ are the Pochhammer symbols and $0 < a,b,c,d \in \R$. These polynomials obey the orthogonality relation
\begin{align}\label{wilson-inner-prod}
\int\limits_0^\infty W_n(x^2\rvert a,b,c,d) W_m(x^2\rvert a,b,c,d) \mathrm{d}\omega(x\rvert a,b,c,d) = N_n(a,b,c,d) \delta_{n,m}.
\end{align}
The weight $\omega(x\rvert a,b,c,d)$ and normalization $N_n(a,b,c,d)$ are given explicitly in \cite{koekoek2010}. For any admissible set of parameters, the Wilson polynomials form a basis of the space of polynomials on the support of $\omega(x\rvert a,b,c,d)$. Belonging to the Askey--Wilson scheme, they are bispectral, that is, they diagonalize a three-term recurrence operator acting on the degree and a difference operator acting on the variable.

The Wilson polynomials form an infinite set of orthogonal polynomials that can be truncated \cite{koekoek2010} to a finite one by setting the parameters as follows
\begin{align*}
 a = \frac{1}{2}(\gamma+\delta+1), \qquad b = \frac{1}{2}(2\alpha-\gamma-\delta+1), \qquad c = \frac{1}{2}(2\beta-\gamma+\delta+1), \qquad d = \frac{1}{2}(\gamma-\delta+1),
\end{align*}
and imposing any of the conditions
\begin{align*}
 \alpha+1=-N, \qquad \beta+\delta+1=-N, \quad \text{or} \quad \gamma+1=-N.
\end{align*}
One thus obtains the Racah polynomials after taking
\begin{align*}
 ix\longmapsto x+\frac{1}{2}(\gamma+\delta+1).
\end{align*}
An additional truncation can be obtained \cite{lemay2016} by imposing
\begin{align}\label{para-racah-truncation}
 a+b+c+d=-N+1.
\end{align}
Indeed, while one is at first sight led to singular expressions, well-defined orthogonal polynomials can nonetheless be obtained through the use of limits and the resulting polynomials, first introduced in \cite{lemay2016}, are the para-Racah polynomials. These polynomials form a three-parameter set of orthogonal polynomials $P_n(x^2\rvert a,c,w)$ of maximal degree $N$. Explicit expressions can be found by setting $N=2j+p$, where $j\in\mathbb{N}$ and $p=0,1$, depending on the parity of $N$. The para-Racah polynomial $P_n(x^2\rvert a,c,w)$ obtained from the truncation \eqref{para-racah-truncation} of the Wilson polynomial $W_n(x^2\rvert a,b,c,d)$ is given by
\begin{align}
 P_n(x^2\rvert a,c,w)=\eta_n \sum\limits_{k=0}^{n}A_{n,k}\varphi_k(x^2), \qquad\qquad \varphi_k(x^2)\equiv (a-ix)_k(a+ix)_k,
\end{align}
where
\begin{align}
 A_{n,k}=
 \begin{cases}
  \frac{(-n)_k(n-N)_k}{(1)_k(-j)_k(a+c)_k(a-c-j+1-p)_k} &k\leq j\\
  \frac{w^{-1}(-n)_k(n-N)_{N-n}(1)_{n+k-1-N}}{(1)_k(-j)_j(1)_{k-j-1}(a+c)_k(a-c-j+1-p)_k} &k>j\\
  0 & k>n,
 \end{cases}
\end{align}
with the normalization given by
\begin{align}
 \eta_n=
 \begin{cases}
  \frac{(1)_n(-j)_n(a+c)_n(a-c-j+1-p)_n}{(-n)_n(n-N)_n} & n\leq j,\\
  \frac{w(1)_n(-j)_j(1)_{n-j-1}(a+c)_n(a-c-j+1-p)_n}{(-n)_n(n-N)_{N-n}(1)_{2n-1-N}} &n>j.
 \end{cases}
\end{align}
These polynomials are orthogonal on a discrete measure that has support on the zeros of the characteristic polynomial $P_{N+1}(x^2\rvert a,c,w)$. The corresponding lattice is a quadratic bi-lattice given by
\begin{align}
 x_{2s+t}=
 \begin{cases}
  -(s+a)^2 \qquad& t=0,\,s=0,1,\dots,j,\\
  -(s+c)^2 \qquad& t=1,\,s=0,1,\dots,j-1+p, 
 \end{cases}
\end{align}
so that
\begin{align}
 \sum\limits_{s=0}^N P_{n}(x^2\rvert a,c,w)P_{m}(x^2\rvert a,c,w)\bar{\omega}_s \propto \delta_{n,m},
\end{align}
where the weight $\bar{\omega}_s$ is given explicitly in \cite{lemay2016}. They also satisfy a three-term recurrence relation and a difference equation. However, they do not appear in classifications of classical orthogonal polynomials as their spectrum is doubly-degenerate.

\subsection{Outline}
The remainder of the paper is organized as follows. In section \ref{section-skop}, the S--Heun operators are introduced and some of their properties are derived. The connection is made with the algebraic Heun operator of the Wilson/Racah type. Section \ref{section-stab} focuses on a subset of the S--Heun operators that preserves the degree of polynomials. A stabilizing algebra is defined from the quadratic relations they obey and its representations are constructed. This algebra is extended to a star algebra in section \ref{section-star} for which a universal presentation is obtained; it is subsequently recognized as a Sklyanin-type algebra. Finally, section \ref{section-rational-sk} provides a presentation of the rational degenerate Sklyanin algebra introduced in \cite{smirnov2010} and gives an isomorphism with the universal algebra of section \ref{section-star}. Using this isomorphism, representations of the rational degenerate Sklyanin algebra on the Wilson and para-Racah polynomials are constructed. A brief conclusion follows.

\section{Sklyanin--Heun operators on a quadratic grid}\label{section-skop}
The generic algebraic Heun operators on a domain $\lambda$ have the property that, when acting on polynomials over $\lambda$, they raise the degree by at most one. The S--Heun operators are a specialization of these Heun operators without a diagonal term. In this section, we first identify the S--Heun operators on the quadratic grid and then proceed with a brief characterization.

\subsection{Sklyanin--Heun operators}
Let $\lambda=\lambda_x$ be a discrete grid indexed by $x$ and define the shift operators $T^{\pm}$ acting on functions on $\lambda$ as follows
\begin{align*}
T^{\pm} f(\lambda_x) \equiv f(\lambda_{x\pm1}).
\end{align*}
Consider a second order operator $S$ with no diagonal term
\begin{align}\label{w-nodiag}
S = A_{1}(\lambda_x) T^{+}+A_{2}(\lambda_x) T^{-},
\end{align}
where $A_1$ and $A_2$ are functions on $\lambda_x$. Demand that $S$ satisfies the degree-raising property
\begin{align}\label{heun-property}
S\cdot p_{n}(\lambda_x) = q_{n+1}(\lambda_x),
\end{align}
for $p_{n}$ and $q_{n+1}$ arbitrary polynomials of degree $n$ and $n+1$, respectively. One can determine the coefficients $A_{1}$ and $A_{2}$ by acting on the first two monomials in $\lambda_x$ as follows
\begin{align}\label{first2monomials}
S\cdot 1 = u_{0}+u_{1}\lambda_x, \quad S \cdot \lambda_x = u_{2}+u_{3}\lambda_x +u_{4} \lambda_x^{2}.
\end{align}
One finds
\begin{align}\label{w-coef}
A_{1}(\lambda_x) &=\frac{u_{2}+u_{3}\lambda_x+u_{4}\lambda_x^{2}-u_{0}\lambda_{x-1}-u_{1}\lambda_{x-1}\lambda_x}{\lambda_{x+1}-\lambda_{x-1}},\\
A_{2}(\lambda_x) &=-\frac{u_{2}+u_{3}\lambda_x+u_{4}\lambda_x^{2}-u_{0}\lambda_{x+1}-u_{1}\lambda_{x+1}\lambda_{x}}{\lambda_{x+1}-\lambda_{x-1}}\nonumber.
\end{align}
The S--Heun operators are defined as the set of operators of the form \eqref{w-nodiag} with the coefficients given in \eqref{w-coef}. As these coefficients admit five independent parameters, the S--Heun operators form a five-dimensional vector space $\mathcal{SH}$ of operators on $\lambda$. A basis for this space can be chosen as follows
\begin{align}
L &= \mathcal{N}(\lambda_x) (T^{+}-T^{-}),\nonumber\\
M_{1} &= \mathcal{N}(\lambda_x) \left[(\lambda_x-\lambda_{x-1})T^{+}+(\lambda_{x+1}-\lambda_{x})T^{-} \right],\nonumber\\
M_{2} &= \mathcal{N}(\lambda_x) \left[(\lambda_x +\lambda_{x-1})T^{+}-(\lambda_{x+1}+\lambda_x)T^{-} \right],\label{sheun-gen}\\
R_{1} &= \mathcal{N}(\lambda_x)\lambda_x \left[(\lambda_x -\lambda_{x-1})T^{+}+(\lambda_{x+1}-\lambda_x)T^{-} \right],\nonumber\\
R_{2} &= \mathcal{N}(\lambda_x)\lambda_x \left[(\lambda_x +\lambda_{x-1})T^{+}-(\lambda_{x+1}+\lambda_x)T^{-} \right],\nonumber
\end{align}
where
\begin{align*}
\mathcal{N}(\lambda_x) \equiv \left[ \lambda_{x+1}-\lambda_{x-1} \right]^{-1}.
\end{align*}
The naming conventions used in \eqref{sheun-gen} will be explained in the next subsection.
\begin{remark}\label{remark-first-order}
 Acting on the left with $T^{+}$ for each of the operators in \eqref{sheun-gen}, it can be seen that the set of operators $\mathcal{SH}$ can also be understood as the set of first order shift operators of step two over $\lambda$ that satisfies the property \eqref{heun-property}.
\end{remark}

\subsection{Sufficiency of the construction}
As established above, for an operator $S$ of the form \eqref{w-nodiag} to satisfy the property \eqref{heun-property}, the expressions \eqref{w-coef} are necessary conditions. The sufficiency of these conditions follows from the ensuing proposition.
\begin{proposition}\label{sufficiency}
 A generic element $S\in \mathcal{SH}$ satisfies the property \eqref{heun-property} if the grid $\lambda_x$ is of one of the following forms
 \begin{align}\label{recurrence-grid}
  \lambda_{x}=\alpha q^{x} + \beta q^{-x} + \kappa, \qquad \lambda_{x}=\alpha x^2 + \beta x + \kappa, \qquad\text{or}\qquad \lambda_{x}=(-1)^x(\alpha x + \beta) + \kappa,
 \end{align}
 for some constants $\alpha$, $\beta$, $\kappa$.
 \begin{proof}
  An element $S\in\mathcal{SH}$ of the form \eqref{w-nodiag} is specified by a set of parameters $\{u_i\}_{i=0,1,\cdots,4}$ \eqref{w-coef}. The action of $S$ on a monomial in $\lambda_x$
  can be reduced by linearity to the five cases given by $u_i=\delta_{i,j}$ for $j=0,1,\dots,4$. Upon inspecting \eqref{w-coef}, one understands that only the operators defined by $u_i=\delta_{i,0}$ or $u_i=\delta_{i,2}$ need to be analyzed; those remaining amount to one of these two operators multiplied by some power of $\lambda_x$.

The first case we treat is $u_i=\delta_{i,2}$ and it corresponds to the operator we have denoted $L$. It can be seen from \eqref{w-coef} that when $u_i=\delta_{i,3}$ or $u_i=\delta_{i,4}$, the corresponding operator is $\lambda_x{}L$ or $\lambda_x{}^2L$, respectively. It follows that for $S$ to satisfy property \eqref{heun-property}, one must have that $L$ decreases the degree of polynomials in $\lambda_x$ by one. Similarly, it follows from \eqref{w-coef} that the case $u_i=\delta_{i,1}$ will satisfy \eqref{heun-property} if the case of $u_i=\delta_{i,0}$, corresponding to the S--Heun operator $\frac{1}{2}(M_1-M_2)$, is an operator that stabilizes the set of polynomials of a given degree.
 
Thus, a generic element of the form \eqref{w-nodiag} will satisfy \eqref{heun-property} if the subset of operators generated by the cases $u_i=\delta_{i,j}$ for $j=0,2,3$ preserves the degree of polynomials. As the generators of this subset are all tridiagonal operators, the proposition follows from the results in \cite{vinet2008} which identifies \eqref{recurrence-grid} as the possible grids allowing second-order difference equations diagonalized by polynomials. 
 \end{proof}
\end{proposition}
On the quadratic grid, it can be shown that a generic element of the vector space spanned by \eqref{sheun-gen} satisfies the property \eqref{heun-property}. Indeed, as derived above, the expressions \eqref{w-coef} for the coefficients are necessary conditions. 
The derivations so far were grid-independent, but to proceed further, one needs to fix the grid. Let us consider the quadratic grid
\begin{align}\label{gridp}
 \lambda_x = x^2.
\end{align}
For this choice of grid, one has that proposition \ref{sufficiency} holds and the sufficiency of the construction is established.

The leading terms of the actions on monomials in $\lambda_x$ are now computed for future reference. In the case of $L$, one obtains
\begin{align}\label{leadterm-L}
L\cdot\lambda_x{}^n = \sum\limits_{k=1}^n \bin{n}{k} \lambda_x^{n-k} \sum\limits_{\substack{j\text{ odd}\\0\leq j\leq k}} \bin{k}{j} (4\,\lambda_x + p^2)^{\frac{j-1}{2}}=n\,\lambda_x{}^{n-1} + O(\lambda_x^{n-2}),
\end{align}
which is verified to be a degree lowering operator. Moreover,
one finds that
\begin{align}
 \frac{1}{2}(M_1-M_2)\cdot \lambda_x{}^n &= \sum\limits_{\substack{k=0\\0\leq j \leq k}}^n \bin{n}{k}\bin{k}{j} \lambda_x^{n-k} \left[\frac{1+(-1)^j}{2}(4\lambda_x + p^2)^{\frac{j}{2}}-\frac{1-(-1)^j}{2}(\lambda_x + 1)(4\lambda_x + p^2)^{\frac{j-1}{2}} \right],\nonumber\\
 &=(1-n)\lambda_x{}^n + O(\lambda_x{}^{n-1}),\label{leadterm-m}
\end{align}
which preserves the degree of polynomials. The actions of the other generators follow from \eqref{leadterm-L} and \eqref{leadterm-m} by noting that
\begin{align}\label{other-sheun}
 \frac{1}{2}(M_1+M_2) = \lambda_x\, L, \qquad R_1 = \lambda_x M_1, \qquad R_2 = \lambda_x M_2.
\end{align}
With the above observations, it follows that a generic linear combination of the basis elements \eqref{sheun-gen} displays the degree raising property \eqref{heun-property}. These calculations enable one to see that the choice of basis \eqref{sheun-gen} decomposes the generic special Heun operator into operators that have a prescribed action on polynomials in $\lambda_x$. Indeed, $L$ can be identified as a lowering operator, $M_{1}$ and $M_{2}$ as stabilizing operators while $R_{1}$ and $R_{2}$ are raising operators. 

\subsection{S--Heun operators of the Wilson type and the Heun--Racah operator}
As the S--Heun operators are specialized algebraic Heun operators \cite{grunbaum2018}, they are related to the general algebraic Heun operators associated to the same grid. The Heun--Racah operator $W$ on the quadratic grid introduced in \cite{bergeron2020a} admits a quadratic embedding in the set $\mathcal{SH}$ of S--Heun operators on the quadratic grid. In view of Remark \ref{remark-first-order}, it will come as no surprise that this embedding is obtained by first conjugating the Heun--Racah operator $W$ by a scaling of the grid $\mu:x\to2x$, such that the shift operators in $W$ act with a step of two. One obtains
\begin{align}\label{heun-wilson-op}
 \mu^{-1}\circ W \circ \mu = R_1(a_1 M_1+a_2 M_2+a_3 L)+R_2(a_4 M_1+a_5 M_2+a_6 L)+a_7 LM_2+a_8 M_2{}^2 +a_9 L^2,
\end{align}
where the coefficients $a_i,\, i=1,2,\dots,9$ are given in terms of the parameters $t_0,t_1,u_0,u_1,u_2,v_0,v_1,v_2$ and $v_3$ of $W$ in \cite{bergeron2020a} as
\begin{align*}
 a_1 &= \frac{t_1+u_2}{4}-\frac{v_3}{16}, & a_2 &= -\frac{t_1}{8}+8u_0+u_1-2v_1+\frac{v_3}{16}, & a_3 &=\frac{1}{4}\left(-8t_0-t_1-64u_0-3u_2+16v_1+2v_2\right),\\
 a_4 &= \frac{u_2}{4}-a_2, & a_5 &=\frac{v_3}{16}, \qquad\quad a_6=a_3-2u_1, & a_7 &= 8 u_0, \qquad\quad a_8 =t_0,\\
 a_9 &= -t_0-24u_0+16v_0. &&&&
\end{align*}

The operator $X$ that acts by multiplication by the grid variable $\lambda_x$ can be written as a quadratic expression in terms of the S--Heun generators:
\begin{align}
 X \equiv x^2 &= 
 (R_1+R_2)(M_1-L)-\frac{1}{2}R_1M_2-\frac{1}{2}R_2M_1.\label{x-op}
\end{align}
 
\section{The stabilizing subalgebra $\mathfrak{stab}$}\label{section-stab}
By direct computations from the definitions \eqref{sheun-gen}, it can be seen that the S--Heun generators satisfy homogeneous quadratic relations, with the complete list given in the appendix \ref{quad-rel}. From these relations, it is observed that the subset of stabilizing S--Heun operators generated by $L,M_1$ and $M_2$ closes as a quadratic algebra to be called $\stab$ whose relations are
\begin{align}\label{stab-sheun}
 [L,M_1]=2L^2,\qquad [L,M_2]=\{M_1,L\}, \qquad [M_1,M_2]=\{M_2,L\}-4L^2.
\end{align}
The Casimir element $C$ is given by
\begin{align}\label{stab-sheun-casimir}
 C=M_1{}^2-\{M_2,L\}+3L^2,
\end{align}
and is equal to the identity in the realization \eqref{sheun-gen} in terms of shift operators. It will prove fruitful to examine the stabilizing algebra \eqref{stab-sheun} in this realization. Knowing that it stabilizes polynomials in $\lambda_x$ of a given degree, one may set up an eigenvalue problem on this space.

\subsection{Diagonalization of a generic linear element}
Consider a generic linear combination of the operators $L,M_1,M_2$
\begin{align}\label{p-op-def}
 P(s,t) = u L + v M_1 + w M_2,
\end{align}
parametrized as follows
$$ u=\frac{(1+2s)(1+2t)-1}{4}, \qquad v=\frac{1}{2}(1+s+t), \qquad w=\frac{1}{2},$$
with $0 < s,t \in \R$ being arbitrary parameters. It is straightforward to show that, under the invertible transformation
$$\rho:x\mapsto -ix,$$
the operator $P$ is given by
$$\tilde{P} \equiv \rho\circ P \circ\rho^{-1} = -\frac{1}{4 i x}\left[(t-ix)(s-ix)\tilde{T}^{+} -(s+ix)(t+ix)\tilde{T}^{-} \right],$$
with $\tilde{T}^{\pm}$ defined by $\tilde{T}^{\pm}f(x) \mapsto f(x\pm i)$. Multiplying each term in the above by $(2 ix \pm 1)/(2 ix \pm 1)$, one recognizes the off-diagonal terms of the difference operator diagonalized by the continuous dual Hahn polynomials \cite{koekoek2010}. Denoting these polynomials as $S_n(x^2 \rvert 1/2,s,t)$ one has
$$ \tilde{P}\, S_n(x^2 \rvert 1/2,s,t) = (n - (s+t)/2)\, S_n(x^2 \rvert 1/2,s,t).$$
Once an element is specified by \eqref{p-op-def}, this defines an eigenbasis in terms of the continuous dual Hahn polynomials. However, no meaningful action can be identified for the remaining elements in $\mathfrak{stab}$. We consider instead quadratic combinations in the elements of the algebra.

\subsection{Action on Wilson polynomials}

A natural action of the stabilizing algebra $\mathfrak{stab}$ on the Wilson polynomials arises from the realization \eqref{sheun-gen}. Indeed, defining the following pair of operators from \eqref{p-op-def}
\begin{align}\label{mu-mustar-def}
\mu^{(a,b,c,d)} =  P(2a-1,2b-1), \qquad\qquad \mu^{*(a,b,c,d)} =  P(2c,2d),
\end{align}
such that manifestly
$$ \mu^{*(a,b,c,d)} = \mu^{(c+1/2,d+1/2,a-1/2,b-1/2)},$$
one has the following proposition: 
\begin{proposition}\label{Qprop}
 The quadratic element $Q\in\mathfrak{stab}$ defined by
 \begin{align}\label{q-op-def}
 Q \equiv \mu^{*(a,b,c,d)}\mu^{(a,b,c,d)},
\end{align}
where $\mu$ and $\mu^*$ are given by \eqref{mu-mustar-def} and with $0 < a,b \in \R$ and $1/2 < c,d \in \R $ is realized, up to a constant term, by the Wilson operator conjugated by the grid scaling
\begin{align}\label{scaling}
  \phi: x\mapsto -2ix. 
\end{align}
\begin{proof}
In the realization \eqref{sheun-gen}, conjugating $Q$ by the scaling transformation \eqref{scaling}, it can be seen by direct calculations that the transformed operator $\tilde{Q}$ is given by
\begin{multline*}
\tilde{Q} \equiv \phi\circ Q \circ \phi^{-1} = \frac{(a-ix)(b-ix)(c-ix)(d-ix)}{2ix(2ix-1)}\,\tilde{T}^{+} + \frac{(a+ix)(b+ix)(c+ix)(d+ix)}{2ix(2ix+1)}\,\tilde{T}^{-}\\
- \left[ \frac{(a-ix)(b-ix)(c-ix)(d-ix)}{2ix(2ix-1)} + \frac{(a+ix)(b+ix)(c+ix)(d+ix)}{2ix(2ix+1)}\right]+(a+b)(c+d-1).
\end{multline*}
The above operator is identified as the Wilson operator \cite{koekoek2010}, up to a constant term.
\end{proof}
\end{proposition}
\begin{remark}
The operator $X$ that acts by multiplication by the variable $\lambda_x$ can be embedded \eqref{x-op} in the set $\mathcal{SH}$ of S--Heun operators. In addition, with the operator $Q$ identified as the Wilson operator, the bispectral pair of operators that generates the Racah/Wilson algebra \cite{granovskii1988,geronimo2010,genest2014} admits an embedding in the set $\mathcal{SH}$ of S--Heun operators. Moreover, a quartic embedding of the Heun--Racah operator \eqref{heun-wilson-op} is obtained from the construction of the Heun--Racah operator \cite{bergeron2020a} by the tridiagonalization \cite{grunbaum2017} of the Racah operator.
\end{remark}
The definition of $Q$ in \eqref{q-op-def} naturally provides a factorization of the Wilson operator in terms of $\mu^{*(a,b,c,d)}$ and $\mu^{(a,b,c,d)}$. Moreover, it directly follows from proposition \ref{Qprop} that the operator $\tilde{Q}$ is diagonalized by the Wilson polynomials:
$$ \tilde{Q}\,W_n(x^2\rvert a,b,c,d) = \left[n(n+a+b+c+d-1)+(c+d)(a+b-1) \right]\,W_n(x^2\rvert a,b,c,d).$$
Introducing a third operator $\tau^{(a,b,c,d)}$ defined by
\begin{align}\label{tau-def}
 \tau^{(a,b,c,d)} = 4L,
\end{align}
a presentation of $\mathfrak{stab}$ in terms of the generators $\mu^{(a,b,c,d)}$, $\mu^{*(a,b,c,d)}$ and $\tau^{(a,b,c,d)}$ can be given for generic values of the parameters $a,b,c,d$. This allows to construct representations of $\mathfrak{stab}$ on the Wilson polynomials.

\begin{proposition}
 A representation of $\mathfrak{stab}$ on the Wilson polynomials $\tilde{W}$ (see \eqref{scaled-wilson}) is given by the following actions
\begin{align}
\mu^{(a,b,c,d)}\cdot\, &\tilde{W}_n(x^2\rvert a,b,c,d) =-(n+a+b-1)\, \tilde{W}_{n}\left(x^2\middle\rvert a-1/2,b-1/2,c+1/2,d+1/2\right),\nonumber\\
 \tau^{(a,b,c,d)}\cdot\, &\tilde{W}_n(x^2\rvert a,b,c,d) =n(n+a+b+c+d-1)\, \tilde{W}_{n-1}\left(x^2\middle\rvert a+1/2,b+1/2,c+1/2,d+1/2\right),\nonumber\\
 \mu^{*\left(a,b,c,d\right)}\cdot\, &\tilde{W}_n(x^2\rvert a,b,c,d) = -\sigma(n+a+b-1)\, \tilde{W}_{n}\left(x^2\middle\rvert a-1/2,b-1/2,c+1/2,d+1/2\right)\label{action-wilson}\\
 &\qquad -(1-\sigma)(n+c+d-1)\, \tilde{W}_{n}\left(x^2\middle\rvert a+1/2,b+1/2,c-1/2,d-1/2\right)\nonumber \\
 &\qquad\qquad +\left[\sigma(ab-cd)-\frac{1}{2}(c+d)-\frac{1}{4}\right]n(n+e_1-1)\, \tilde{W}_{n-1}\left(x^2\middle\rvert a+1/2,b+1/2,c+1/2,d+1/2\right),\nonumber\\
 &\nonumber\\
 \sigma\equiv (a&+b-c-d)^{-1}, \qquad e_1 \equiv a+b+c+d, \nonumber
\end{align}
with
\begin{align}\label{scaled-wilson}
\tilde{W}_n(x^2\rvert a,b,c,d) \equiv \phi^{-1}\cdot W_n(x^2\rvert a,b,c,d) = W_n\Big(-\frac{x^2}{4}\Big\rvert a,b,c,d\Big),
\end{align}
and $\phi$ defined in \eqref{scaling}.
\begin{proof}
 The conjugation of the three operators $\mu^{(a,b,c,d)}$, $\mu^{*(a,b,c,d)}$ and $\tau^{(a,b,c,d)}$ by the scaling map \eqref{scaling} yields operators that are identified as the structure and forward shift operators for the Wilson polynomials \cite{koekoek2010}. These structure operators have a known action on the Wilson polynomials \cite{miller1987}. Using the identity
 \begin{align}\label{mustar-decomp}
   \mu^{*(a,b,c,d)}=\sigma\,\mu^{(a,b,c,d)}+(1-\sigma)\mu^{(c,d,a,b)}+\left[\sigma(ab-cd)-\frac{1}{2}(c+d)-\frac{1}{4}\right]\tau^{(a,b,c,d)},
 \end{align}
which is directly verified and applying the scaling \eqref{scaling} to the polynomials to get \eqref{scaled-wilson}, one obtains the actions \eqref{action-wilson}. As one can use the orthogonality relation \eqref{wilson-inner-prod} of the Wilson polynomials to express all polynomials with shifted parameters in \eqref{action-wilson} as sums of Wilson polynomials with the initial parameters, these actions define representations of the stabilizing algebra $\mathfrak{stab}$ on the Wilson polynomials.
\end{proof}
\end{proposition}
\section{Extension of $\mathfrak{stab}$ to a star algebra}\label{section-star}
The construction laid out in the preceding section parallels the structural approach to orthogonal polynomials due to Kalnins and Miller \cite{miller1987,kalnins1989}. In particular, Miller derives in \cite{miller1987} the orthogonality \eqref{wilson-inner-prod} of the Wilson polynomials from the structural recurrence relations associated to $\mu^{(a,b,c,d)}$ and $\tau^{(a,b,c,d)}$ by identifying the operator $\mu^{*(a,b,c,d)}$ and deriving an inner product such that this operator is the adjoint of $\mu^{(a,b,c,d)}$. An operator $\tau^{*(a,b,c,d)}$ is then identified as the adjoint of $\tau^{(a,b,c,d)}$. A similar approach in the context of the S--Heun operators can be pursued at the algebraic level.

The representations defined through \eqref{action-wilson} are endowed with a natural inner product inherited from the orthogonality relation \eqref{wilson-inner-prod}. This enables one to define a star operation, such that $\mu^{*(a,b,c,d)}$ is precisely the adjoint of $\mu^{(a,b,c,d)}$ under the inner product. It follows that $\tilde{Q}$ is a self-adjoint operator. However, the stabilizing algebra is not closed under the star operation. This can be seen by taking the adjoint of $\tau^{(a,b,c,d)}$, a lowering operator, which would involve raising operators that are not contained in the stabilizing algebra $\mathfrak{stab}$. We shall now extend $\mathfrak{stab}$ to its closure under the star operation.

\subsection{Star operation}
With the help of \eqref{wilson-inner-prod}, one constructs an operator as the adjoint of the forward shift operator. This leads to the backward shift operator for the Wilson polynomials \cite{koekoek2010} with action given by
\begin{align}\label{taustar-action}
 (\phi^{-1}\circ\tau^{*(a,b,c,d)}\circ\phi)\cdot\,W_n(x^2\rvert a,b,c,d) = W_{n+1}\left(x^2\middle\rvert a-1/2,b-1/2,c-1/2,d-1/2\right).
\end{align}
The operator $\tau^{*(a,b,c,d)}$ can then be decomposed in terms of the S--Heun operators as follows
\begin{align}\label{taustar-def}
 \tau^{*(a,b,c,d)} = a_1 L + a_2 M_1 + a_3 M_2 + a_4 R_1 + a_5 R_2,
\end{align}
with the coefficients given by
\begin{align*}
 a_1 &= 4e_4-e_3+\frac{e_1-1}{4}, & a_2 &= e_3-\frac{e_2}{2}+\frac{e_1}{8}, & a_3 &=\frac{e_2}{2}-\frac{5e_1}{8}+\frac{1}{2}, & a_4 &=\frac{e_1}{4}-\frac{3}{8}, & a_5 = -\frac{1}{8},
\end{align*}
where $e_1$, $e_2$, $e_3$ and $e_4$ are the elementary symmetric polynomials in the four parameters $a,b,c$ and $d$:
\begin{align}\label{elementary-pol}
\begin{split}
 e_1 &= a+b+c+d, \qquad\qquad\quad e_2 = ab+ac+ad+bc+bd+cd,\\
 e_3 &= abc+abd+acd+bcd, \,\quad e_4 = abcd.
\end{split}
\end{align}

Introducing $\tau^{*(a,b,c,d)}$ as a fourth generator together with those of the stabilizing algebra $\mathfrak{stab}$ leads to an algebra closed under the star operation.
\begin{proposition}
 The algebra $\struc$ generated by $\mu^{(a,b,c,d)}$, $\mu^{*(a,b,c,d)}$, $\tau^{(a,b,c,d)}$ and $\tau^{*(a,b,c,d)}$, together with the relations induced from their definitions in terms of S--Heun operators given in \eqref{mu-mustar-def}, \eqref{tau-def} and \eqref{taustar-def} admits the natural star map defined from its canonical action on the generators:
\begin{align}\label{star-map}
 *\,:\,\, &\tau\longmapsto\tau^*,\\
  &\mu\longmapsto\mu^*.
\end{align}
\begin{proof}
 The result follows from the results of \cite{miller1987} after conjugation of the generators by the scaling map \eqref{scaling}. 
\end{proof}
\end{proposition}

\subsection{A universal presentation of $\mathfrak{stab}^*$}
The algebra $\struc$ can be presented in terms of quadratic relations by making use of the relations given in the appendix \ref{quad-rel}. However, such a presentation obfuscates the structure of the algebra because the parameters $a,b,c$ and $d$ of the Wilson polynomial appear explicitly in the relations. Thus, it does not define uniquely an algebra associated to the quadratic grid.

Recall that the normalized Wilson polynomials are known \cite{miller1987} to be fully symmetric under permutations of their four parameters. However, the definitions for the two stabilizing generators given in \eqref{mu-mustar-def} do not make this symmetry manifest, because they contain the specific parameters of the representation. Nervertheless, the permutation symmetry of the polynomials can be made manifest at the level of the algebra to obtain a universal presentation.
\begin{proposition}\label{univ-pres}
 The algebra $\mathfrak{stab}^*$ admits a presentation as a unital associative algebra with four generators $U,V,Y$ and $R$ obeying the following relations
 \begin{align}
 [V,Y]&=-\{U,Y\}, \qquad [U,Y]=-\{Y,Y\}, \qquad [U,V]=\{V,Y\}-2\{Y,Y\},\nonumber\\
 [R,Y]&=\{U,U\}-\{U,V\}+\{V,Y\}, \qquad\qquad\,\, [R,V]=2\{V,Y\}-\{Y,Y\}-\{V,V\}-\{U,R\},\label{star-alg-def}\\
 [R,U]&=\{U,V\}+2\{V,Y\}-2\{U,Y\}-\{V,V\}-\{Y,Y\}-\{R,Y\}.\nonumber
\end{align}
The two Casimir operators are given by
\begin{align}\label{star-alg-cas}
 Q_1 = U^2-\{V,Y\}+3Y^2, \qquad Q_2 = U^2 + V^2 - \{U,V\} - \{U,Y\} - \{R,Y\}.
\end{align}
\begin{proof}
Consider the following generic linear combination of generators
$$ u\,\mu^{(a,b,c,d)} + v\,\mu^{*(a,b,c,d)}.$$
Acting with the symmetric group $S_4$ on the parameters $(a,b,c,d)$, one constructs a fully symmetric element in terms of the S--Heun operators as follows
\begin{multline*}
 \frac{1}{|S_4|}\sum\limits_{\sigma \in S_4} \left[u\,\mu^{\sigma(a,b,c,d)} + v\,\mu^{*\,\sigma(a,b,c,d)}\right]=
 \\ \frac{1}{2}[(u-v) - e_1(u+v)]M_1 -\frac{1}{2}(u+v)M_2+\left[\frac{e_1}{2}(u-v)-\frac{2e_2}{3}(u+v)\right]L.
\end{multline*}
Setting $u=1$ and either $u=v$ or $u=-v$ in the above yields two independent generators that are manifestly symmetric and can be used instead of $\mu$ and $\mu^*$ to obtain another presentation of $\mathfrak{stab}^*$. The relations in this new presentation now only involve the elementary symmetric polynomials \eqref{elementary-pol}. Subsequently, it becomes straightforward to eliminate all remaining parameters in the algebraic relations by further redefining the generators as
\begin{align}\label{uvy-def}
 U = M_1+e_1\,L, \qquad V = M_2 + e_1\,M_1 + \frac{1}{2}e_1{}^2L, \qquad Y=L,
 \end{align}
 \begin{multline}\label{r-def}
 R = R_2 + (2 e_1-3)R_1 + \frac{1}{2}(3e_1{}^2-10e_1+4)M_2 + \frac{1}{2}(e_1+1)(e_1{}^2-4e_1+2)M_1\\
 +\frac{1}{8}(e_1{}^4 - 4e_1{}^3-8e_1{}^2+24e_1-8)L.
\end{multline}
Using the quadratic relations of the S--Heun operators given in the appendix \ref{quad-rel}, the relations \eqref{star-alg-def}, as well as the centrality of the two operators in \eqref{star-alg-cas}, are verified.
\end{proof}
\end{proposition}
In a realization in terms of S--Heun operators, the Casimir operators \eqref{star-alg-cas} are proportionnal to the identity and the coefficients are functions of the parameters of the polynomials. One has
\begin{align}
 Q_1 = 1, \qquad Q_2 = (e_1-2)(e_1-4),
\end{align}
where $e_1$ is given in \eqref{elementary-pol}.
\begin{remark}
 While a universal presentation of $\mathfrak{stab}^*$ has been given in proposition \ref{univ-pres}, the star structure is not universal and depends explicitly on the representation parameters. This is not surprising because the map \eqref{star-map} is constructed using the inner product \eqref{wilson-inner-prod} corresponding to a specific realization with fixed parameters. Nevertheless, one can work in a specific realization and write the generators in \eqref{star-alg-def} in terms of the structural operators \eqref{mu-mustar-def}, \eqref{tau-def} and \eqref{taustar-def} as follows
\begin{align}\label{staralg-on-struc}
\begin{split}
 U &= \frac{1}{|S_4|}\sum\limits_{\sigma \in S_4} \left[\mu^{\sigma(a,b,c,d)}-\mu^{*\,\sigma(a,b,c,d)}\right], \qquad V = \frac{1}{|S_4|}\sum\limits_{\sigma \in S_4} \left[-\mu^{\sigma(a,b,c,d)}-\mu^{*\,\sigma(a,b,c,d)}\right]+\alpha\, Y,\\
 Y &= \frac{1}{4}\tau, \qquad R = 8 \tau^{*(e_1,e_2,e_3,e_4)} - (2-3\alpha)V+(1-3\alpha+\beta)U + (1-\beta +\gamma) Y,
\end{split}
\end{align}
where
\begin{align}\label{alpha-beta-gamma}
 \alpha = \frac{1}{2}e_1{}^2-\frac{4}{3}e_2, \qquad \beta = -e_1{}^3+4e_2 e_1-8e_3, \quad\text{and}\quad \gamma = \frac{3}{4}\alpha\,e_1{}^2 - e_1{}^2 e_2 + 8e_1e_3-32e_4,
\end{align}
with $e_1,e_2,e_3$ and $e_4$ given in \eqref{elementary-pol}. With the above, one obtains
\begin{align*}
\begin{split}
 U^*&=-U, \qquad Y^*=\frac{1}{32}\left[R + (2-3\alpha)V-(1-3\alpha+\beta)U - (1-\beta +\gamma)Y\right], \qquad V^* = V+\alpha\,(Y^*-Y),\\
 R^* &= \big[32+\alpha\,(2-3\alpha)\big] Y - (2-3\alpha)V - (1-3\alpha+\beta)U + \big[1-\beta +\gamma-\alpha\,(2-3\alpha)\big] Y^*.
\end{split}
\end{align*}
\end{remark}

\subsection{The algebra $\mathfrak{stab}^{*}$ as a Sklyanin algebra}
It can be seen from \eqref{uvy-def} and \eqref{r-def} that the generators of $\mathfrak{stab}^*$ only depend on the parameters $a,b,c,d$ via the elementary symmetric polynomial $e_1(a,b,c,d)$. Thus, they will be invariant under a commensurate increase and decrease of any pair of parameters. A glance at \eqref{taustar-def} indicates that this will not be the case for $\tau^{*(a,b,c,d)}$. However, a pseudo-commutation relation similar to the one introduced by Rains in \cite{rains2006} is obtained.
\begin{proposition}
 In the realization \eqref{sheun-gen} the identity
 \begin{align}\label{rains-relation}
 \tau^{*(a,b,c+k,d-k)}\, \tau^{*(a+\frac{1}{2},b+\frac{1}{2},c-\frac{1}{2},d-\frac{1}{2} )}=\tau^{*(a,b,c,d)}\,\tau^{*(a+\frac{1}{2},b+\frac{1}{2},c-\frac{1}{2}+k,d-\frac{1}{2}-k)},
\end{align}
 is satisfied. Moreover, at the abstract level \eqref{rains-relation} encodes the algebraic relations of the $\mathfrak{stab}^*$ algebra \eqref{star-alg-def}.
 \begin{proof}
  Using the definition \eqref{taustar-def}, the identity \eqref{rains-relation} is readily verified. The second statement is demonstrated by using \eqref{staralg-on-struc} to express $\tau^{*(a,b,c,d)}$ in terms of the generators \eqref{uvy-def} and \eqref{r-def} as
\begin{align*}
8\tau^{*(a,b,c,d)} = R + (2-3\alpha)V-(1-3\alpha+\beta)U - (1-\beta +\gamma) Y,
\end{align*}
where $\alpha,\beta$ and $\gamma$ are given in \eqref{alpha-beta-gamma}. Upon using the above in \eqref{rains-relation}, one can pick any one of the parameters $a,b,c,d$ and take the remaining ones to be vanishing. Equating the coefficients of each power of the remaining non-zero parameter in the left- and right-hand side of \eqref{rains-relation} yields a set of relations that is algebraically identical to the relations \eqref{star-alg-def}.
 \end{proof}

\end{proposition}
That the relations of $\ustruc$ are encoded in the identity \eqref{rains-relation} identifies the $\ustruc$ algebra as a Sklyanin-type algebra \cite{rains2006}.

\section{The rational degenerate Sklyanin algebra}\label{section-rational-sk}
The rational degenerate Sklyanin algebra $\skr$ is obtained in \cite{smirnov2010} from the Sklyanin algebra \cite{sklyanin1984} and is associated to a rational degeneration of an elliptic $R$-matrix. A presentation can be given as a unital associative algebra generated by four elements $S_0,S_3,S_+,S_-$ obeying the defining relations
\begin{align}\label{rdskly-def}
\begin{split}
 [S_0,S_-] &= -2\{S_-,S_-\}, \qquad [S_0,S_+] = 16\{S_3,S_-\}-16\{S_-,S_-\}+2\{S_+,S_-\}- 4 \{S_3,S_3\},\\
 [S_+,S_-] &= 2\{S_0,S_3\}, \qquad\quad\,\,\, [S_0,S_3] = 2\{S_3,S_-\}-8\{S_-,S_-\}, \qquad [S_3,S_\pm] = \pm\{S_0,S_\pm\}.
\end{split}
\end{align}
The rational degenerate Sklyanin algebra admits two Casimir operators which are given in the above presentation by
\begin{align}\label{rdskly-casimir}
 C_1 = S_0^2 + S_3^2 + \frac{1}{2} \{S_+,S_-\}, \qquad C_2 = \frac{1}{2}\{S_+,S_-\}+2\{S_-,S_3\}+S_3^2-6\{S_-,S_-\}.
\end{align}
The presentation \eqref{rdskly-def} is recovered from the one in \cite{smirnov2010} upon setting the free parameter $\eta=1$ and defining $S_\pm = S_1\pm iS_2$. The following proposition identifies the $\ustruc$ algebra with the rational degenerate Sklyanin algebra.
\begin{proposition}\label{iso-sklyanin}
 The $\skr$ algebra defined in \eqref{rdskly-def} is isomorphic to the $\ustruc$ algebra defined in \eqref{star-alg-def}.
 \begin{proof}
  The following map is readily verified to be an isomorphism of algebras.
  \begin{align}\label{sklyanin-on-staralg}
 S_0=4Y-4U, \qquad S_3 = 4U-2Y-4V, \qquad S_+ = 16R-14Y-8U+24V, \qquad S_- = -2Y.
\end{align}
 \end{proof}
\end{proposition}
\subsection{A realization in terms of difference operators}
A realization of the rational degenerate Sklyanin algebra in terms of difference operators is provided in \cite{smirnov2010}. The Casimir elements are realized as multiples of the identity and are given by
$$ C_1 = 16(2s+1)^2  Id, \qquad C_2=64 s(s+1) Id.$$
The generators thus represented can be written in terms of the S--Heun operators \eqref{sheun-gen} as follows
$$ S_0 = 4(2s-1)L -4 M_1, \qquad S_3 = -2(2s-1)^2 L + 4(2s-1) M_1 - 4 M_2, \qquad S_1 - i S_2 = -2 L, $$
$$ S_1 + i S_2 = -2(4s^2-1)(4s^2-8s-1)L -8(2s-1)(4s^2-4s-1)M_1 + 8(2s-1)(6s+1)M_2 - 16(4s-1)R_1 + 16 R_2. $$
It is immediate from the above that the realization in terms of S--Heun operators of the $\skr$ algebra involves coefficients that depend on the values of the Casimir operators. A similar observation could be made for the case of the $\ustruc$ algebra in \eqref{uvy-def} and \eqref{r-def}. It follows from proposition \ref{iso-sklyanin} that the parameters $e_1$ and $s$ are related by
\begin{align*}
 e_1 = 2-2s.
\end{align*}

\subsection{A family of representations}
The identification of the rational degenerate Sklyanin algebra $\skr$ with the $\ustruc$ algebra directly leads to a family of representations of $\skr$ on the Wilson polynomials.
\begin{proposition}\label{sklyanin-rep}
 A representation of the rational degenerate Sklyanin algebra $\skr$ \eqref{rdskly-def} on the  Wilson polynomials is defined by the following actions
 \begin{multline*}
S_0 \cdot \tilde{W}_{n}\left(x^2\middle\rvert a,b,c,d\right) =4\sigma (n+c+d-1)\, \tilde{W}_{n}\left(x^2\middle\rvert a+1/2,b+1/2,c-/2,d-1/2\right)\\
+(4\sigma(ab - cd)-e_1)n(n+e_1-1)\,\tilde{W}_{n-1}\left(x^2\middle\rvert a+1/2,b+1/2,c+1/2,d+1/2\right)\\
-4\sigma (n+a+b-1)\, \tilde{W}_{n}\left(x^2\middle\rvert a-1/2,b-1/2,c+1/2,d+1/2\right),
\end{multline*}
\begin{multline*}
 S_3\cdot \tilde{W}_{n}\left(x^2\middle\rvert a,b,c,d\right) = -4(n+a+b-1)\,\tilde{W}_{n}\left(x^2\middle\rvert a-1/2,b-1/2,c+1/2,d+1/2\right)\\
 +\frac{1}{2}(8(ab+cd)-e_1{}^2-1)n(n+e_1-1)\,\tilde{W}_{n-1}\left(x^2\middle\rvert a+1/2,b+1/2,c+1/2,d+1/2\right)\\
 -4(n+c+d-1)\,\tilde{W}_{n}\left(x^2\middle\rvert a+1/2,b+1/2,c-/2,d-1/2\right),\\
\end{multline*}
\begin{multline*}
 S_- \cdot \tilde{W}_{n}\left(x^2\middle\rvert a,b,c,d\right) = -\frac{1}{2}n(n+e_1-1)\,\tilde{W}_{n-1}\left(x^2\middle\rvert a+1/2,b+1/2,c+1/2,d+1/2\right),\shoveleft{}
\end{multline*}
\begin{multline*}
 S_+ \cdot \tilde{W}_{n}\left(x^2\middle\rvert a,b,c,d\right) = 128\tilde{W}_{n+1}\left(x^2\middle\rvert a-1/2,b-1/2,c-1/2,d-1/2\right)\\
 +8(6\alpha-1+2\beta\sigma)(n+a+b-1)\,\tilde{W}_{n}\left(x^2\middle\rvert a-1/2,b-1/2,c+1/2,d+1/2\right)\\
 +8(6\alpha-1-2\beta\sigma)(n+c+d-1)\,\tilde{W}_{n}\left(x^2\middle\rvert a+1/2,b+1/2,c-/2,d-1/2\right)\\
 +8\left[(1-6\alpha)(ab+cd)-2\beta\sigma(ab-cd)+\xi\right]n(n+e_1-1)\,\tilde{W}_{n-1}\left(x^2\middle\rvert a+1/2,b+1/2,c+1/2,d+1/2\right),
\end{multline*}
where $\alpha,\beta$ and $\gamma$ are defined in \eqref{alpha-beta-gamma} and
\begin{align*}
 \xi \equiv \frac{1}{2}(1-2e_1{}^2+e_1{}^4-256e_4),
\end{align*}
with $e_1$ and $e_4$ defined in \eqref{elementary-pol}.

\begin{proof}
One first derives the action of the symmetrized structure operators on the Wilson polynomials. It can be seen from \eqref{action-wilson} and \eqref{taustar-action} that the expressions in the case of $\tau$ and $\tau^*$ are fully symmetric under permutations of the parameters such that their actions are invariant under the symmetrization. To obtain similar expressions for $\mu$ and $\mu^*$, one uses \eqref{mustar-decomp} to write
\begin{multline}\label{mu-mustar-delta}
 \mu^{(a,b,c,d)}\pm\mu^{*(a,b,c,d)} =(\mu^{(a,b,c,d)}\pm\mu^{(c,d,a,b)})\mp\frac{1}{2}\big(c+d+\frac{1}{2}\big)\tau^{(a,b,c,d)} \\
 \pm \sigma\left[\mu^{(a,b,c,d)}-\mu^{(c,d,a,b)}+(ab-cd)\tau^{(a,b,c,d)} \right].
\end{multline}
The last term in the right-hand side of \eqref{mu-mustar-delta} is independent of the parameters as
\begin{align}\label{pam-free-part}
 \sigma\left[\mu^{(a,b,c,d)}-\mu^{(c,d,a,b)}+(ab-cd)\tau^{(a,b,c,d)} \right] = Y-U,
\end{align}
and is thus invariant under the symmetrization. As it is verified that
\begin{align*}
 \frac{1}{\rvert S_4 \rvert} \sum\limits_{\pi\in S_4}(\mu^{\pi(a,b,c,d)}-\mu^{\pi(c,d,a,b)})=0,
\end{align*}
one can use the invariance of $\tau$ under permutations of the parameters to obtain from \eqref{mu-mustar-delta} using \eqref{pam-free-part} that
\begin{align*}
  \frac{1}{\rvert S_4 \rvert} \sum\limits_{\pi\in S_4}(\mu^{\pi(a,b,c,d)}-\mu^{*\,\pi(a,b,c,d)}) = \sigma\mu^{(c,d,a,b)} - \sigma\mu^{(a,b,c,d)} + \left[\frac{1}{4}(e_1 +1) -\sigma(ab-cd) \right]\tau^{(a,b,c,d)}.
\end{align*}
Likewise, observing that
$
 \mu^{(a,b,c,d)}+\mu^{(c,d,a,b)}+(ab+cd)\tau^{(a,b,c,d)},
$
is symmetric under permutations of the parameters, one can use the invariance of $\tau$ and \eqref{pam-free-part} in \eqref{mu-mustar-delta} to obtain
\begin{multline*}
 \frac{1}{\rvert S_4 \rvert} \sum\limits_{\pi\in S_4}(\mu^{\pi(a,b,c,d)}+\mu^{*\,\pi(a,b,c,d)})=(1+\sigma)\mu^{(a,b,c,d)}+(1-\sigma)\mu^{(c,d,a,b)}\\
 +\left[(ab+cd)+\sigma(ab-cd) - \frac{1}{3}e_2 -\frac{1}{4}(e_1+1) \right]\tau^{(a,b,c,d)}.
\end{multline*}
The actions on the scaled Wilson polynomials \eqref{scaled-wilson} of $\tau,\tau^*$ and of the operators in \eqref{mu-mustar-delta} are obtained from \eqref{action-wilson} and found to be:
\begin{multline*}
 \frac{1}{\rvert S_4 \rvert} \sum\limits_{\pi\in S_4}(\mu^{\pi(a,b,c,d)}-\mu^{*\,\pi(a,b,c,d)})\cdot\, \tilde{W}_n(x^2\rvert a,b,c,d) =\\
  \left[\frac{1}{4}(e_1 +1) -\sigma(ab-cd) \right]n(n+e_1-1) \tilde{W}_{n-1}(x^2\rvert a+1/2,b+1/2,c+1/2,d+1/2)\\
 +\sigma(n+a+b-1)\tilde{W}_n(x^2\rvert a-1/2,b-1/2,c+1/2,d+1/2)\\
 -\sigma(n+c+d-1)\tilde{W}_n(x^2\rvert a+1/2,b+1/2,c-1/2,d-1/2),
\end{multline*}
\begin{multline*}
 \frac{1}{\rvert S_4 \rvert} \sum\limits_{\pi\in S_4}(\mu^{\pi(a,b,c,d)}+\mu^{*\,\pi(a,b,c,d)})\cdot\, \tilde{W}_n(x^2\rvert a,b,c,d) =\\
 \left[(ab+cd)+\sigma(ab-cd) - \frac{1}{3}e_2 -\frac{1}{4}(e_1+1) \right]n(n+e_1-1) \tilde{W}_{n-1}(x^2\rvert a+1/2,b+1/2,c+1/2,d+1/2)\\
 +(\sigma-1)(n+c+d-1)\tilde{W}_n(x^2\rvert a+1/2,b+1/2,c-1/2,d-1/2)\\
 -(\sigma+1)(n+a+b-1)\tilde{W}_n(x^2\rvert a-1/2,b-1/2,c+1/2,d+1/2)
\end{multline*}
\begin{align*}
\frac{1}{\rvert S_4 \rvert} &\sum\limits_{\sigma\in S_4}\tau^{\sigma(a,b,c,d)}\cdot\, \tilde{W}_n(x^2\rvert a,b,c,d) =n(n+e_1-1)\, \tilde{W}_{n-1}\left(x^2\middle\rvert a+1/2,b+1/2,c+1/2,d+1/2\right),&\\
 \frac{1}{\rvert S_4 \rvert} &\sum\limits_{\sigma\in S_4}\tau^{*\,\sigma(a,b,c,d)}\cdot\,\tilde{W}_n(x^2\rvert a,b,c,d) = \tilde{W}_{n+1}\left(x^2\middle\rvert a-1/2,b-1/2,c-1/2,d-1/2\right).&
\end{align*}
Using \eqref{staralg-on-struc} and the above, one can construct a representation of \eqref{star-alg-def} on the Wilson polynomials \eqref{scaled-wilson}. Proposition \ref{sklyanin-rep} then follows from proposition \ref{iso-sklyanin}.
\end{proof}
\end{proposition}

Finite-dimensional representations can be obtained by truncating the representations of proposition \ref{sklyanin-rep}.
\begin{proposition}\label{para-racah-rep}
 The finite-dimensional representations obtained from truncations of the representations in proposition \ref{sklyanin-rep} act on the para-Racah polynomials.
 \begin{proof}
Looking at the content of proposition \ref{sklyanin-rep}, it is seen that the only generator that raises the degree is $S_+$. Using \eqref{sklyanin-on-staralg}, \eqref{staralg-on-struc} and \eqref{r-def} this degree-raising action can be traced back to the following combination of S--Heun operators
\begin{align*}
 R_2+(2e_1-3)R_1.
\end{align*}
With the help of \eqref{other-sheun},\eqref{leadterm-m} and \eqref{leadterm-L}, one can obtain the leading term of the action of the above operator on a polynomial of degree $N$ in $\lambda_x$:
\begin{align}
& R_1 \cdot \lambda_x{}^N = \lambda_x{}^{N+1}+O(\lambda_x{}^{N}),  \qquad R_2 \cdot \lambda_x{}^N = (2N-1)\lambda_x{}^{N+1}+O(\lambda_x{}^N),\nonumber\\
 & \left[R_2+(2e_1-3)R_1\right]\cdot \lambda_x{}^N =2(N-1+e_1)\lambda_x{}^{N+1}+O(\lambda_x{}^N).
\end{align}
Demanding that the leading term in the above vanishes is tantamount to truncating the representation of proposition \ref{sklyanin-rep} at the degree $N$. This truncation condition is precisely the one that leads to the para-Racah polynomials \eqref{para-racah-truncation}. Thus, the finite-dimensional representations of the rational degenerate Sklyanin algebra obtained under this truncation have for basis the para-Racah polynomials.
 \end{proof}
\end{proposition}
The actions of the generators in these truncated representations are as given in proposition \ref{sklyanin-rep}, although one has to carry the appropriate limiting process described in \cite{lemay2016} to deal with the otherwise singular expressions. Proposition \ref{para-racah-rep} provides for the algebraic interpretation of the para-Racah polynomials as the basis elements of the finite-dimensional representations of the rational degenerate Sklyanin algebra.

\section{Conclusion}
This paper has introduced the S--Heun operators associated to the quadratic grid as a special case of the algebraic Heun operator. These operators were shown to form a five-dimensional space. The subset of these operators which stabilizes the space of polynomials of a given degree was identified and the algebra that they realize was examined. The extension of this stabilizing algebra to a star algebra was identified as the rational degenerate Sklyanin algebra. This definition of the rational degenerate Sklyanin algebra through S--Heun operators directly led to the construction of infinite-dimensional representations on the Wilson polynomials as well as finite-dimensional representations on the para-Racah polynomials.

The rational degenerate Sklyanin algebra is known \cite{smirnov2010} to be a one parameter deformation of the Yangian $Y(\mathfrak{sl}_2)$. In the same way that the Yangian $Y(\mathfrak{sl}_2)$ is the quantum algebra that encodes the symmetry of integrable $XXX$ spin-half chains associated with the ordinary rational $R$-matrix, the rational degenerate Sklyanin algebra can be understood as the symmetry algebra of a generalized $XXX$ chain corresponding to a deformed rational $R$-matrix, a new integrable model. Thus, it would be of interest to use the representations introduced in section \ref{section-rational-sk} to construct explicit realizations of this new integrable model in terms of finite and infinite spin chains. In the finite case, one would expect the para-Racah polynomials to appear as the basis of representations of the symmetry algebra. Interestingly, these para polynomials were first introduced in the context of perfect state transfer on spin chains \cite{vinet2012} and the advances in this paper suggest they would also find applications as solutions to new integrable spin chain models.

\section*{Acknowledgments}
The authors would like to thank Jean-Michel Lemay for useful discussions. While part of this research was conducted, GB held a scholarship from the Institut des Sciences Mathématiques of Montreal and JG held an Alexander-Graham-Bell scholarship from the Natural Sciences and Engineering Research Council of Canada (NSERC). The research of LV is funded in part by a Discovery Grant from NSERC. AZ gratefully holds a CRM-Simons professorship and his work is supported by the National Science Foundation of China (Grant No.11771015).


\section{Appendix}

\subsection{Quadratic relations of the S--Heun operators}\label{quad-rel}
The set of homogeneous quadratic algebraic relations satisfied by the S--Heun operators is given below for reference:

$$[L,M_1]=2L^2,\qquad [L,M_2]=\{M_1,L\}, \qquad [M_1,M_2]=\{M_2,L\}-4L^2, $$
$$ [L,R_1]=M_1^2+L^2+\{M_1,L\}+\frac{1}{2}\{M_2,L\}, \qquad [L,R_2]=M_1^2+L^2+\{M_1,L\}+\frac{1}{2}\{M_2,L\}+\{M_1,M_2\},$$
$$ [M_1,R_1]=2M_1^2-3L^2+\{M_1,M_2\}-\frac{1}{2}\{M_1,L\}-\{M_2,L\},$$
$$[M_1,R_2]=M_1^2 +M_2^2+7L^2 +2\{R_2,L\}-\frac{5}{2}\{M_1,L\}-5\{M_2,L\},$$
$$[R_1,M_2]=3L^2-M_1^2-M_2^2+2\{R_1+R_2,L\}-\{R_1,M_1\}-\{M_1,M_2\}-5\{M_1,L\}-\frac{9}{2}\{M_2,L\}, $$
$$ [R_2,M_2]=Y^2-M_1^2-M_2^2+\{R_1,M_1-M_2\}-\{M_1,M_2\}+\frac{1}{2}\{M_1,L\}, $$
$$
 [R_2,R_1]=2R_1^2+M_1^2+2M_2^2+3L^2+\frac{1}{2}\{R_2-R_1,L\}-\frac{3}{2}\{R_1+R_2,M_2\} + \{M_1,M_2\}+\frac{3}{2}\{M_1,L\}-\frac{1}{2}\{M_2,Y\},$$
 
 $$ M_1^2-\{M_1,M_2\}+3L^2=1, \qquad \{R_1-R_2,L\}+M_2^2+\{M_2,L\}-3L^2=-3,$$
 $$ -2\{R_1,L\}-3L^2+\{M_1,M_2\}+2\{M_1+M_2,L\}=4,$$
 $$ M_1^2+\frac{1}{2}L^2+\{R_1,M_1-M_2\}-\frac{5}{2}\{R_1,L\}-2\{R_2,L\}+\{R_1+R_2,M_1\}+\frac{1}{4}\{M_1,M_2\}+6\{M_1,L\}+4\{M_2,L\}=0.$$

\newpage
\bibliographystyle{abbrv}
\bibliography{s-heun-quad.bib}

\end{document}